\documentclass[times]{article}

\usepackage{moreverb}

\usepackage[dvips,colorlinks,bookmarksopen,bookmarksnumbered,citecolor=red,urlcolor=red]{hyperref}
\usepackage{comment}
\usepackage{graphicx}
\usepackage{amsmath,amsfonts,amssymb,amsbsy,empheq,stmaryrd}
\usepackage{MnSymbol}
\usepackage{subfigure}
\usepackage{multirow}
\usepackage{tabularx}
\usepackage{xcolor}
\usepackage[utf8]{inputenc}
\usepackage{tikz}
\usetikzlibrary{plotmarks,mindmap,trees,shapes,arrows,shapes.multipart}

\usepackage[a4paper,body={18cm,26cm},top=2cm]{geometry}

\makeatletter
\newif\if@restonecol
\makeatother

\usepackage[algo2e,ruled,lined,titlenumbered]{algorithm2e}

\DeclareMathAlphabet{\mathdj}{U}{msb}{m}{n}
\DeclareMathAlphabet{\mathbf}{OML}{cmm}{b}{it}
\DeclareMathAlphabet{\mathbfsl}{OT1}{cmr}{bx}{sl}

\newcommand{\range}{\ensuremath{\operatorname{Range}}}
\renewcommand{\ker}{\ensuremath{\operatorname{Ker}}}
\newcommand{\R}{\ensuremath{\mathbb{R}}}

\newcommand{\sU}{\ensuremath{\mathcal{U}}}

\newcommand{\s}{\ensuremath{^{(s)}}}
\newcommand{\sT}{\ensuremath{^{(s)^T}}}

\newcommand{\ddd}{\ensuremath{^{\diamondbackslash}}}
\newcommand{\ddl}{\ensuremath{^{\diamondminus}}}
\newcommand{\ddc}{\ensuremath{^{\diamondvert}}}
\newcommand{\dddT}{\ensuremath{^{\diamondbackslash^T}}}
\newcommand{\ddlT}{\ensuremath{^{\diamondminus^T}}}

\newcommand{\eftrace}{\ensuremath{\mathbf{t}}}

\newcommand{\intp}{\ensuremath{\Gamma_{\!A}}}
\newcommand{\intd}{\ensuremath{\Gamma_{\!B}}}

\newcommand{\efassemp}{\ensuremath{\mathbf{A}}}
\newcommand{\efassemd}{\ensuremath{\mathbf{B}}}

\newcommand{\dep}{\ensuremath{u}}

\newcommand{\depv}{\ensuremath{v}}
\newcommand{\res}{\ensuremath{r}}
\newcommand{\efdep}{\ensuremath{\mathbf{\dep}}}

\newcommand{\efres}{\ensuremath{\mathbf{\res}}}

\newcommand{\efdepi}{\ensuremath{\mathbf{\dep}_b}}
\newcommand{\ovefdepi}{\ensuremath{\mathbf{\dep}_{\!A}}}
\newcommand{\efdepik}{\ensuremath{\mathbf{\dep}_{b_k}}}
\newcommand{\ovefdepik}{\ensuremath{\mathbf{\dep}_{\!A_k}}}
\newcommand{\efddepik}{\ensuremath{\mathbf{\overset{\circ}{\dep}}_{b_k}}}
\newcommand{\efddepk}{\ensuremath{\mathbf{\overset{\circ}{\dep}}_{k}}}

\newcommand{\ovefddepik}{\ensuremath{\mathbf{\overset{\circ}{\dep}}_{\!A_k}}}
\newcommand{\efdepv}{\ensuremath{\mathbf{\depv}}}
\newcommand{\efshape}{\ensuremath{\mathbf{N}}}
\newcommand{\eps}{\ensuremath{\varepsilon}}
\newcommand{\efl}{\ensuremath{\boldsymbol{\lambda}}}
\newcommand{\eflam}{\ensuremath{\boldsymbol{\lambda}_b}}
\newcommand{\uveflam}{\ensuremath{\boldsymbol{\lambda}_{\!B}}}

\newcommand{\eflamk}{\ensuremath{\boldsymbol{\lambda}_{b_k}}}
\newcommand{\uveflamk}{\ensuremath{\boldsymbol{\lambda}_{B_k}}}
\newcommand{\efdlamk}{\ensuremath{\boldsymbol{\overset{\circ}{\lambda}}_{b_k}}}
\newcommand{\uvefdlamk}{\ensuremath{\boldsymbol{\overset{\circ}{\lambda}}_{\!B_k}}}

\newcommand{\uveflamz}{\ensuremath{\boldsymbol{\lambda}_{\!B_0}}}

\newcommand{\efmu}{\ensuremath{\boldsymbol{\mu}_b}}
\newcommand{\efmuk}{\ensuremath{\boldsymbol{\mu}_{b_k}}}
\newcommand{\efdmuk}{\ensuremath{\boldsymbol{\overset{\circ}{\mu}}_{b_k}}}
\newcommand{\impe}{\ensuremath{\mathbf{Q}_b}}
\newcommand{\efxi}{\ensuremath{\mathbf{x}_b}}

\newcommand{\efrhs}{\ensuremath{\mathbf{b}}}
\newcommand{\efx}{\ensuremath{\mathbf{x}}}
\newcommand{\proj}{\ensuremath{\mathbf{P}_B}}

\newcommand{\forcf}{\ensuremath{f}}
\newcommand{\forcF}{\ensuremath{F}}
\newcommand{\efforce}{\ensuremath{\mathbf{\forcf}}}
\newcommand{\bdf}{\ensuremath{\partial_{\forcf}\Omega}}
\newcommand{\bdu}{\ensuremath{\partial_{\dep}\Omega}}

\newcommand{\efstiff}{\ensuremath{\mathbf{K}}}
\newcommand{\efker}{\ensuremath{\mathbf{R}}}
\newcommand{\efalp}{\ensuremath{\boldsymbol{\alpha}}}

\newcommand{\ovefdepii}{\ensuremath{\mathbf{v}_{\!A}}}

\newcommand{\ovefdepiik}{\ensuremath{\mathbf{v}_{\!A_k}}}

\newcommand{\uveflamm}{\ensuremath{\boldsymbol{\gamma}_{\!B}}}
\newcommand{\ovefddepiik}{\ensuremath{\mathbf{\overset{\circ}
{v}}_{\!A_k}}}
\newcommand{\uvefdlammk}{\ensuremath{\boldsymbol{\overset{\circ}{\gamma}}_{B_k}}}

\newcommand{\stiff}{\ensuremath{\mathbf{K}}}
\newcommand{\schurp}{\ensuremath{\mathbf{S}}}
\newcommand{\schurd}{\ensuremath{\mathbf{F}}}
\newcommand{\schurm}{\ensuremath{\mathbf{M}}}

\newcommand{\matzero}{\ensuremath{\mathbfsl{0}}}

\title{Substructured formulations of nonlinear structure problems --  Influence of the interface condition}

\author{Camille Negrello$^{1}$, Pierre Gosselet$^{1}$, Christian Rey$^{1,2}$ and Julien Pebrel$^{1}$ \\
{(1)} LMT Cachan, ENS Cachan/CNRS/Univ. Paris Saclay, \\
61 Avenue du Pr\'esident Wilson, 94235 Cachan France.\\
{(2)} Safran Tech, rue des jeunes bois, Chateaufort,\\
CS 80112, 78772 Magny-les-Hameaux France.
}

\begin{document}

\maketitle

\begin{abstract}
We investigate the use of non-overlapping domain decomposition (DD) methods for nonlinear structure problems. The classic techniques would  combine a global Newton solver with a  linear DD solver for the tangent systems. We propose a framework where we can swap Newton and DD, so that we solve independent nonlinear problems for each substructure and linear condensed interface problems. The objective is to decrease the number of communications between subdomains and to improve parallelism. Depending on the interface condition, we derive several formulations which are not equivalent, contrarily to the linear case. Primal, dual and mixed variants are described and assessed on a simple plasticity problem.

\end{abstract}

\section{Introduction}
In order to solve large nonlinear structure problems, an efficient  strategy is to combine a Newton-based solver which leads to a sequence of linear systems, and a domain decomposition approach to solve the tangent systems. Indeed such a strategy combines well-known and robust methods for which many refinement are available: Newton can be tangent/constant/secant/modified/arc-length  \cite{Kelley03}, and domain decomposition solvers  \cite{Man93,Far94bis,Let94,FARHAT:2001:FETI_DP,GOSSELET.2007.1} can be equipped with preconditioners and coarse problems which make them reliable and scalable \cite{KLAWONN:2007:SCAL_FETIDP,KLAWONN:2008:JAGGED,spillane2013abstract,GOSSELET.2015.1.1}. Moreover, the computation of stiffness matrices being done in parallel independently for each subdomain, and information being reusable from one system to another \cite{Rey98,Ris00,GOSSELET.2012.1}, the overall performance is very satisfying in general \cite{Far00,BRANDS:2008:BIO_FETI_DP}.

But there are cases where such an approach is not as pertinent as expected. For instance, when dealing with strong localized nonlinearity, many (global)  Newton iterations are required, whereas most of the structure undergoes a linear evolution. In this situation, it would be interesting to differentiate between the linear or nonlinear nature of each subdomain's behavior, in order to decrease the number of global iterations and communications. The possibility to conduct local nonlinear computations has been investigated for a long time in the Schwarz framework (with or without overlap between subdomains) \cite{badea91,Dryja97,Lad99d,cai02,Lad07,Hwang07}. More recently the possibility to define nonlinear versions of the Schur complement methods (often called nonlinear relocalization techniques) was studied: BDD in \cite{Cresta07,BORDEU.2009.3.1}, FETI in \cite{Peb08}, FETI-DP in \cite{klawonn14}, BDDC in \cite{hinojosa14,klawonn14}.

The aim of this paper is to give a formal framework to develop the nonlinear versions of the well-known linear solver FETI \cite{Far94bis}, BDD \cite{Man93,Let94} and FETI2LM\cite{FXFETI2LM09}. They rely on the concept of nonlinear condensation which was exposed in conferences \cite{PEBREL.2009.3.2} and seminars \cite{GOSSELET.2011.6.x}.

In section~\ref{sec:ref} and~\ref{sec:ssf}, we present the nonlinear system in its monolithic and substructured form. In section~\ref{sec:nlcond}, we introduce the concept of nonlinear condensation in its primal, dual and mixed versions. In section~\ref{sec:solstrat}, we show that applying a Newton method to these condensed formulations leads, at each outer iteration, to the parallel solution of local nonlinear systems with, depending on the formulation, Dirichlet, Neumann or Robin boundary conditions, and to the linear interface systems which can be solved by classical BDD, FETI or FETI2LM. Iterative solvers being involved inside the outer loop, we interpret the method as an inexact Newton solver \cite{dembo1982inexact} which leads us, in section~\ref{sec:error}, to tune the convergence thresholds depending on the current residual in order to avoid oversolving. First assessments are given in section~\ref{section:assess} on a simple problem, yet representative of the nonlinearity encountered in industrial problems; a discussion ends the paper.

\section{Reference problem}\label{sec:ref}
We consider the classic problem of the evolution, under the small perturbation hypothesis, of a structure occupying the domain $\Omega$, submitted to body forces $\forcf$, to traction forces $\forcF$ on the part $\bdf$ of its boundary and to given displacements $\dep_g$ on the complementary part $\bdu\neq\emptyset$. Note that the small perturbation hypothesis is crucial for the handling of rigid body motions in the dual approach, but it can be relaxed in other cases.

Let $\eps(\dep)$ be the symmetric part of the gradient of displacement and $\sigma$ the Cauchy stress tensor.
The problem to be solved can be written as:
\begin{equation*}
\begin{aligned}
&\text{at time }t\in[0,T],\ \text{find }\dep\in \sU\ /\ \forall \depv\in \sU_0 \\
&\int_\Omega \sigma:\eps(\depv) d\Omega= \int_\Omega \depv^T\forcf d\Omega+ \int_{\bdf} \depv^T\forcF dS \\
&\sigma = \sigma(\varepsilon(\dep(\tau)),\tau\in [0,t])
\end{aligned}
\end{equation*}
where $\sU$ is the space of kinematically admissible fields and $\sU_0$ is the associated vector space:
\begin{equation*}
\sU=\left\{\dep\in H^1(\Omega),\ \dep =\dep_g \ \text{on}\ \bdu
\right\}
\end{equation*}
The notation $\sigma(\varepsilon(\dep(\tau)),\tau\in [0,t])$ means that the stress at one point depends on the whole history of the strain at this point (local nonlinear behavior). This history is most often materialized by internal variables like anelastic strain, hardening or damage.

The problem is discretized in space using the finite element method. The domain $\Omega$ is meshed; let $\efshape$ be the matrix of shape functions, and $\efdep$ be the nodal displacement unknowns such that $\dep=\efshape\efdep$. The problem is also supposed to be discretized in time, the discrete reference problem to be solved at Step $t_n$ can be written as:
\begin{equation}\label{eq:form}\text{Find }\efdep(t_n)\text{ so that }
\efforce_{int}(\efdep)+\efforce_{ext}=0
\end{equation}
with
\begin{equation*}
\begin{aligned}
&\efforce_{ext}=\int_\Omega \efshape^T\forcf d\Omega+ \int_{\bdf} \efshape^T\forcF dS \\
&\efdepv^T\efforce_{int}=-\int_\Omega\sigma_h:\eps(\efshape\efdepv) d\Omega\\
&\sigma_h = \sigma_h(\varepsilon(\efshape\efdep(t_j)),\ j\leqslant n)
\end{aligned}
\end{equation*}
$\sigma_h$ is the discrete counterpart of $\sigma$, it depends of the whole discrete history.  
Integrals are classically computed numerically using Gauss quadrature, so that stress and internal variables are defined at Gauss points. Throughout the rest of the document, the time step will not be mentioned. In order to shorten expressions, Dirichlet boundary conditions $\dep_g$ are implicitly taken into account within $\efforce_{int}$ and $\efforce_{ext}$.

Note that in the linear case we have:  $\efforce_{int}(\efdep)=-\stiff\efdep$ where $\stiff$ is the stiffness matrix.

The classical solution strategy to Problem \eqref{eq:form} is to use a Newton-Raphson algorithm to linearize the problem, and then solve a sequence of tangent systems.
\section{Substructured formulation}\label{sec:ssf}
\begin{figure}[ht]\centering
 \includegraphics[width=0.6\textwidth]{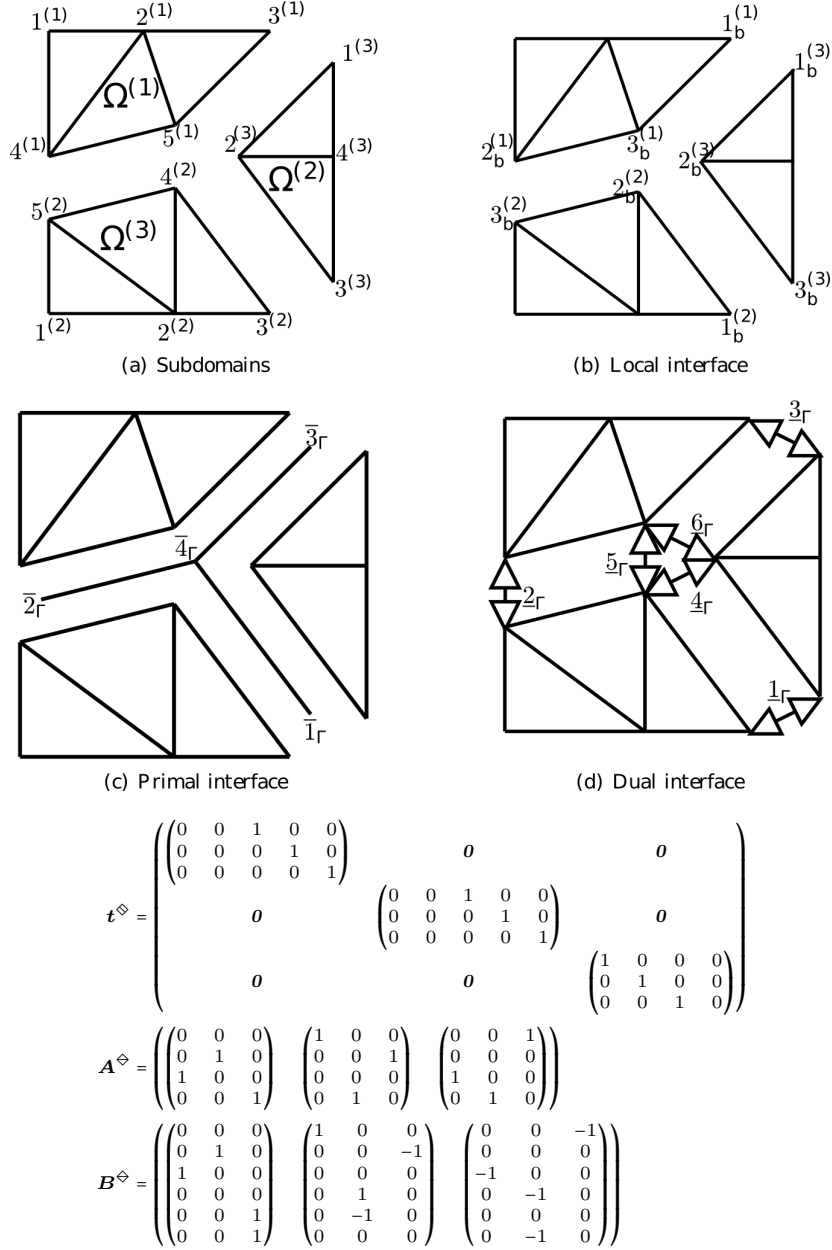}
  \scriptsize\begin{equation*}
 \begin{aligned}
 \eftrace\ddd&=\begin{pmatrix}\begin{pmatrix}0&0&1&0&0\\0&0&0&1&0\\0&0&0&0&1\end{pmatrix}&\matzero&\matzero\\
 \matzero&\begin{pmatrix}0&0&1&0&0\\0&0&0&1&0\\0&0&0&0&1\end{pmatrix}&\matzero\\
 \matzero&\matzero&\begin{pmatrix}1&0&0&0\\0&1&0&0\\0&0&1&0\end{pmatrix}\end{pmatrix}\\
 \efassemp\ddl&=\begin{pmatrix}\begin{pmatrix}0&0&0\\0&1&0\\1&0&0\\0&0&1\end{pmatrix}&\begin{pmatrix}1&0&0\\0&0&1\\0&0&0\\0&1&0\end{pmatrix}&\begin{pmatrix}0&0&1\\0&0&0\\1&0&0\\0&1&0\end{pmatrix}\end{pmatrix}\\
 \efassemd\ddl&=\begin{pmatrix}\begin{pmatrix}0&0&0\\0&1&0\\1&0&0\\0&0&0\\0&0&1\\0&0&1\end{pmatrix}&\begin{pmatrix}1&0&0\\0&0&-1\\0&0&0\\0&1&0\\0&-1&0\\0&0&0\end{pmatrix}&\begin{pmatrix}0&0&-1\\0&0&0\\-1&0&0\\0&-1&0\\0&0&0\\0&-1&0\end{pmatrix}\end{pmatrix}
 \end{aligned}\end{equation*}\normalsize
 \caption{Local numberings, interface numberings, trace and assembly operators}\label{fig:operators}
 \end{figure}

We consider the conforming partition of $\Omega$ into $N$ non-overlapping subdomains $\Omega\s$, so that each element exactly belongs to one subdomain. Superscript $(s)$ will refer to data attached to domain $\Omega\s$. In order to ease the treatment of groups of subdomains, we define the following block notations:
\begin{equation*}\begin{aligned}
\mathbf{x}\ddc&=\begin{pmatrix}\vdots\\\mathbf{x}\s\\\vdots\end{pmatrix},&
\mathbf{x}\ddl&=\begin{pmatrix}\ldots&\mathbf{x}\s&\ldots\end{pmatrix},&
\mathbf{x}\ddd&=\begin{pmatrix}\ddots&&\mathbfsl{0}\\&\mathbf{x}\s&\\\mathbfsl{0}&&\ddots\end{pmatrix}
\end{aligned}\end{equation*}
Note that in the case of nonlinearly-dependent data, we use the same notation but the dependence is implicitly local:
\begin{equation*}
\efforce_{int}\ddc(\efdep\ddc)=\begin{pmatrix}\vdots\\\efforce_{int}\s(\efdep\s)\\\vdots\end{pmatrix}
\end{equation*}
The nodes on the interface between two subdomains play a specific role which we need to highlight.
We note $\Gamma^{(i,j)}$ the set of nodes shared by subdomains $\Omega^{(i)}$ and $\Omega^{(j)}$, and $\Gamma\s=\bigcup_j\Gamma^{(s,j)}$ the interface nodes of subdomain $\Omega\s$. We use the subscript $b$ for nodes belonging to the interface $\Gamma\s$ and the subscript $i$ for internal nodes. The trace operator $\eftrace\s$ extracts interface nodal values (on $\Gamma\s$) from subdomain data (in $\Omega\s$):
\begin{equation*}\begin{aligned}
&\forall s\in\llbracket 1,N\rrbracket,\ \eftrace\s\efdep\s=\efdepi\s \text{ which also can be written as }\ \eftrace\ddd\efdep\ddc=\efdepi\ddc\\
&\text{assuming adapted ordering, }\eftrace\s=\begin{pmatrix} \mathbfsl{0}_{bi}\s & \mathbf{I}_{bb}\s\end{pmatrix}
\end{aligned}
\end{equation*}
See Figure~\ref{fig:operators}(a,b) for an example. We note $\R_b\ddc=\range(\eftrace\ddd)$ the space of local interface vectors to which any vector $\efxi\ddc$ belongs.

Let $\intp=\bigcup\Gamma\s$ denote the totality of interface nodes, and let $\efassemp\s$ be the canonical operator which injects nodes from $\Gamma\s$ to $\intp$. $\efassemp\s$ is a boolean full column-rank matrix. With one local interface node in $\Gamma\s$ associated with exactly one node in $\intp$, the global assembling operator $\efassemp\ddl$ has no left-kernel. It is then a full row-rank matrix (sometimes called a primal assembling operator). See Figure~\ref{fig:operators}(b,c) for an example. We note $\R^{\intp}=\range(\efassemp\ddl)$ the space of vectors defined on that ``primal'' interface.

Any matrix $\efassemd\ddl$ satisfying $\range(\efassemd\ddlT)=\ker(\efassemp\ddl)$ can be used as a dual assembling operator. In Figure~\ref{fig:operators}(b,d) we give the most classical choice where $\efassemd\ddl$ is a signed boolean matrix which describes each connection between interface nodes; we note $\intd$ the set of connections. Note that operator $\efassemd\ddl$ needs not to be full row-rank (in the classical case it is not whenever one interface degree of freedom is shared by more than two subdomains).  We note $\R^{\intd}=\range(\efassemd\ddl)$ the space of vectors defined on that ``dual'' interface.

In order to decouple equations betweens subdomains, we introduce the nodal reaction $\eflam\s$ imposed on subdomain $\Omega\s$ by its neighbors. 
The substructured reference problem can be written as:
\begin{equation}\label{eq:subform}
\begin{aligned}
&\text{Find }(\efdep\ddc,\eflam\ddc)\text{ such that}\\
&\efforce_{int}\ddc(\efdep\ddc)+\efforce\ddc_{ext}+\eftrace\dddT\eflam\ddc=0\\
&\efassemp\ddl \eflam\ddc=0\\
&\efassemd\ddl \efdepi\ddc=0
\end{aligned}
\end{equation}
where the first equation expresses the equilibrium of each subdomain submitted to given forces and to unknown interface reactions $\eflam\ddc$.  The primal assembling operator $\efassemp\ddl$ enables us to express the action-reaction principle while the dual assembling operator $\efassemd\ddl$ enables us to express the continuity of the displacement field. In addition to Figure~\ref{fig:operators}, note that in the case of two subdomains these operators can be written as $\efassemp\ddl \eflam\ddc=\eflam^{(1)}+\eflam^{(2)}$ and $\efassemd\ddl \efdepi\ddc=\efdepi^{(1)}-\efdepi^{(2)}$.

Assembling operators satisfy the following relationships:
\begin{itemize}
\item Assembling operators are orthogonal in the following sense:
\begin{equation}\label{eq:orthoass1}
\efassemp\ddl\efassemd\ddlT=0
\end{equation}
\item Assembling operators generate local interface nodal vectors:
\begin{equation}\label{eq:orthoass1a}
\range(\efassemp\ddlT)\overset{\perp}{\oplus}\range(\efassemd\ddlT) = \R_b\ddc
\end{equation}
\item Any local interface vector is uniquely defined as a combination of a balanced vector $\efassemd\ddlT \mathbf{x}_{\!B} $ and a continuous vector $\efassemp\ddlT \mathbf{x}_{\!A}$:
\begin{equation}\label{eq:orthoass2}
\begin{aligned}
\forall \mathbf{x}_b\ddc&,\ \exists (\mathbf{x}_{\!B},\mathbf{x}_{\!A})\in\R^{\intd}\times\R^{\intp}
/ \ \mathbf{x}_b\ddc = \efassemd\ddlT \mathbf{x}_{\!B} + \efassemp\ddlT \mathbf{x}_{\!A}\\
&\text{indeed } \left\{\begin{aligned}
\mathbf{x}_{\!A}&=\left(\efassemp\ddl\efassemp\ddlT\right)^{-1}\efassemp\ddl\mathbf{x}_b\ddc\\
\mathbf{x}_{\!B}&=\left(\efassemd\ddl\efassemd\ddlT\right)^{+}\efassemd\ddl\mathbf{x}_b\ddc
\end{aligned}\right.
\end{aligned}
\end{equation}
The use of the pseudo-inverse $\left(\efassemd\ddl\efassemd\ddlT\right)^{+}$ is due to the potential presence of redundancies in the description of the connectivity between subdomains, though it is applied to a vector which belongs to $\range{\left(\efassemd\ddl\right)}=\range{\left(\efassemd\ddl\efassemd\ddlT\right)}$ so that $\mathbf{x}_B$ is well defined, and $\efassemd\ddlT \mathbf{x}_{\!B} $ does not depend on the choice of the pseudo-inverse. 
\end{itemize}
\section{Nonlinear condensations}\label{sec:nlcond}
The substructured formulation \eqref{eq:subform} is strictly equivalent to the global formulation \eqref{eq:form}. We now propose various solution strategies which (under certain assumptions) all converge to the reference solution, though these methods, which also employ Newton-Raphson algorithm, are meant to generate a different sequence of linear systems. 

\subsection{Primal formulation}
The primal formulation consists in rewriting system \eqref{eq:subform} in terms of one unknown interface displacement field $\ovefdepi$ :
\begin{equation}\label{eq:primform}
\begin{aligned}
&\text{Find }\ovefdepi\in \R^{\intp}\text{ such that } \efassemp\ddl \eflam\ddc=0 \\
&\text{where } \eflam\ddc:= -\left[ \efforce_{int}\ddc(\efdep\ddc)+\efforce\ddc_{ext} \right]_b \\
&\text{and }\efdep\ddc\text{ solves }\left\{\begin{aligned} &\left[\efforce_{int}\ddc(\efdep\ddc)+\efforce\ddc_{ext}\right]_i=0\\
&\eftrace\ddd \efdep\ddc = \efassemp\ddlT \ovefdepi\end{aligned}\right.\\
\end{aligned}
\end{equation}
The last set of equations corresponds to the solution to independent mechanical problems for each subdomain with imposed displacement at the interface (a Dirichlet condition which implies that the continuity is automatically insured: $\efdepi\ddc=\efassemp\ddlT \ovefdepi\Rightarrow \efassemd\ddl \efdepi\ddc=0$). This displacement has to be found so that the associated reactions $\eflam\ddc$ are balanced on the interface.

If we assume that the last set of equations has a unique solution for any imposed interface displacement, then we can define an operator $\schurp\s_{nl}$ so that:
\begin{equation}
\eflam\s=\schurp\s_{nl}({\efassemp\s}^T\ovefdepi;\efforce_{ext}\s)
\end{equation}
This operator can be viewed as a nonlinear version of the Schur complement; it computes the reaction associated with a given displacement. In the linear case, an explicit expression can be given:
\begin{equation}\label{eq:primal_linear}
\begin{aligned}
\schurp_l\s(\efdepi\s;\efforce_{ext}\s)&= \schurp_t\s\efdepi\s - \efrhs_p\s \\
\text{with }\schurp_t\s&=\stiff\s_{bb}-\stiff\s_{bi}{\stiff\s_{ii}}^{-1}\stiff\s_{ib}\\
\efrhs_p\s&={\efforce_{ext}\s}_b-\stiff\s_{bi}{\stiff\s_{ii}}^{-1}{\efforce_{ext}\s}_{i}
\end{aligned}
\end{equation}

The primal nonlinear condensed system reads:
\begin{empheq}[box=\fbox]{equation}\label{eq:primal_condensed}
\begin{aligned}
&\text{Find }\ovefdepi\in\R^{\intp}
\text{ such that}\\
&\efassemp\ddl\schurp_{nl}\ddc(\efassemp\ddlT\ovefdepi;\efforce_{ext}\ddc)=0
\end{aligned}
\end{empheq}

\subsection{Dual formulation}
The dual formulation consists in rewriting system \eqref{eq:subform} in terms of one unknown interface reaction field:
\begin{equation}\label{eq:dualform}
\begin{aligned}
&\text{Find }\uveflam\in\R^{\intd} \text{ such  that }\efassemd\ddl \eftrace\ddd\efdep\ddc=0 \\
&\text{ where }\efdep\ddc\text{ solves }\efforce_{int}\ddc(\efdep\ddc)+\efforce\ddc_{ext}+\eftrace\dddT\efassemd\ddlT\uveflam=0
\end{aligned}
\end{equation}
The last equation corresponds to the solution to independent mechanical problems for each subdomain with imposed traction at the interface (local Neumann problems). This traction has to be found so that displacements are continuous at the interface. Implicitly local reactions have been defined by $\eflam\ddc=\efassemd\ddlT\uveflam$, which guarantees the interface equilibrium $\efassemp\ddl \eflam\ddc=0$.

In order for problem  \eqref{eq:dualform} to be well-posed, one necessary condition is the self-equilibrium of the substructure. Let $\efker\s$ be the basis of kinematically admissible rigid body motions of subdomain $\Omega\s$ (if they exist, these displacements are exactly represented in the finite element function space), then we have ${\efker\s}^T\efforce_{int}\s=0$ (internal forces do not develop work within rigid body motions). This implies the admissibility condition for $\uveflam$:
\begin{equation}\label{eq:admiss}
\efker\dddT\left(\efforce\ddc_{ext}+\eftrace\dddT\efassemd\ddlT\uveflam\right)=0
\end{equation}

If we assume that \eqref{eq:dualform} has a unique solution (up to a rigid body motion $\efker\ddd\efalp\ddc$, where $\efalp\ddc$ is the unknown amplitude of rigid body motions) for any imposed interface traction, then we can define an operator $\schurd_{nl}\s$ such that:
\begin{equation}
\efdepi\s=\schurd\s_{nl}({\efassemd\s}^T\uveflam;\efforce_{ext}\s)+\eftrace\s\efker\s\efalp\s
\end{equation}
This operator can be viewed as a nonlinear version of the dual Schur complement (as employed in classical FETI methods \cite{Far94} of which we borrow the notation $\schurd$); it computes the displacement associated with  given reaction. In the linear case, an explicit expression can be given:
\begin{equation}\label{eq:dual_linear}
\begin{aligned}
\schurd_l\s(\eflam\s;\efforce_{ext}\s)&= \schurd_t\s\eflam\s + \efrhs_d\s\\ 
\text{with }\schurd_t\s&={\eftrace\s}{\stiff\s}^+{\eftrace\s}^T\\
\efrhs_d\s&={\eftrace\s}{\stiff\s}^+\efforce_{ext}\s
\end{aligned}
\end{equation}
We recall the following classical relationships \cite{GOSSELET.2007.1}:
\begin{equation}
  \begin{aligned}
    \efker\s &= \ker(\efstiff\s) \\ 
    \efker\s_b &= \eftrace\s\efker\s = \ker(\schurp_t\s) \\ 
    {\efker\s} ^T \efforce\s_{ext} &= {\efker\s_b}^T \efrhs_p\s
  \end{aligned}
\end{equation}

Finally, the dual nonlinear condensed system reads:
\begin{empheq}[box=\fbox]{equation}\label{eq:dual_condensed}
\begin{aligned}
&\text{Find }\uveflam\in\R^{\intd}
,\efalp\ddc\text{ such that}\\
&\left\{\begin{aligned}
&\efassemd\ddl\left(\schurd\ddc_{nl}(\efassemd\ddlT\uveflam;\efforce_{ext}\ddc)+\eftrace\ddd\efker\ddd\efalp\ddc\right)=0\\
&\efker\dddT\left(\efforce\ddc_{ext}+\eftrace\dddT\efassemd\ddlT\uveflam\right)=0
\end{aligned}\right.
\end{aligned}
\end{empheq}

\subsection{Mixed formulation}
The mixed formulation consists in introducing a new interface variable for each subdomain $\efmu\s$:
\begin{equation}
  \efmu\ddc=\eflam\ddc+\impe\ddd\efdepi\ddc
\end{equation}
The symmetric positive definite matrix $\impe\ddd$ is a parameter of the method. It can be interpreted as a stiffness (or an impedance) added to the interface (like in Robin-type boundary conditions). 
 
The properties of the matrix $\impe\ddd$ imply a partitioning of the space of subdomain interface vectors:
\begin{equation}\label{eq:mixedortho1}
\begin{aligned}
\ker(\efassemd\ddl)&\oplus\ker(\efassemp\ddl\impe\ddd)=\R_b\ddc
\end{aligned}
\end{equation}
This property enables us to reformulate both interface conditions \eqref{eq:subform} in one single boundary equation:
\begin{equation}
  \efassemp\ddlT\left( \efassemp\ddl \impe\ddd\efassemp\ddlT \right)^{-1} \efassemp\ddl  \efmu\ddc - \efdepi\ddc = 0
\end{equation}
indeed when multiplying this equation on the left by $\efassemd\ddl$ or $\efassemp\ddl\impe\ddd$ we recover classical interface conditions. 

The mixed formulation of Problem \eqref{eq:subform} can be written as:
\begin{equation}\label{eq:mixedform}
\begin{aligned}
&\text{Find }\efmu\ddc\in\R_b\ddc\text{ such that }\efassemp\ddlT\left( \efassemp\ddl \impe\ddd\efassemp\ddlT \right)^{-1} \efassemp\ddl  \efmu\ddc - \efdepi\ddc = 0\\
&\text{where }\efdep\ddc\text{ solves }\efforce_{int}\ddc(\efdep\ddc)-\eftrace\dddT\impe\ddd\eftrace\ddd\efdep\ddc+\eftrace\dddT\efmu\ddc+\efforce\ddc_{ext}=0
\end{aligned}
\end{equation}
The last equation corresponds to the solution of independent nonlinear problems for each subdomain under Robin boundary conditions. If we assume that it has a unique solution for any given $\efmu\ddc$ (the addition of matrix $\impe\ddd$ suppresses rigid body motions), then we can define a nonlinear interface operator:
\begin{equation}
\efdepi\s=\schurm\s_{nl}(\efmu\s;\efforce_{ext}\s,\impe\s)
\end{equation}
In the linear case, the operator can be written as:
\begin{equation}\label{eq:mixed_linear}
\begin{aligned}
\schurm\s_{l}(\efmu\s;&\efforce_{ext}\s,\impe\s)= \schurm\s_{t}\efmu\s + \efrhs_m\s \\
\text{with }\schurm\s_{t} &=\eftrace\s(\stiff\s+{\eftrace\s}^T\impe\s{\eftrace\s})^{-1}{\eftrace\s}^T\\
\efrhs_m\s &= \eftrace\s(\stiff\s+{\eftrace\s}^T\impe\s{\eftrace\s})^{-1}\efforce_{ext}\s
\end{aligned}
\end{equation}
One can recognize the operators of the FETI2LM approach \cite{FXFETI2LM09}.

The mixed nonlinear condensed problem can be written as:
\begin{empheq}[box=\fbox]{equation}\label{eq:mixed_condensed}
\begin{aligned}
\text{Find }\efmu\ddc\in\R_b\ddc\text{ such that}& \\
  \efassemp\ddlT\left( \efassemp\ddl \impe\ddd\efassemp\ddlT \right)^{-1} \efassemp\ddl  \efmu\ddc &- \schurm\ddc_{nl}(\efmu\ddc;\efforce_{ext}\ddc,\impe\ddd) = 0
\end{aligned}
\end{empheq}

Note that the method makes use of the operator $\left( \efassemp\ddl \impe\ddd\efassemp\ddlT \right)$ which has exactly the structure of an assembled condensed stiffness matrix. The factorization of such a matrix is expensive in the general case, but we can choose a specific fill-in of matrix $\impe\ddd$ for the assembled matrix $\left( \efassemp\ddl \impe\ddd\efassemp\ddlT \right)$ to have a block diagonal structure, which makes its handling much cheaper.\medskip

Note that an equivalent formulation is possible where the boundary unknowns $\efmu\ddc$ are replaced by two interface unknowns: one balanced force ${\uveflamm}$ and one continuous displacement fields ${\ovefdepii}$:
\begin{equation}\label{eq:mixedortho3}
  \efmu \ddc = \efassemd\ddlT \uveflamm + \impe\ddd\efassemp\ddlT \ovefdepii
\end{equation}
This is due to \eqref{eq:mixedortho1} which gives the following decomposition, similar to \eqref{eq:orthoass2}:
\begin{equation}\label{eq:mixedortho2}
\begin{aligned}
  \forall \efxi\ddc&\in\R_b\ddc, \exists  (\mathbf{x}_B,\mathbf{x}_A)\in\R^{\intd}\times\R^{\intp}
\ /\ \efxi\ddc = \impe\ddd\efassemp\ddlT \mathbf{x}_A + \efassemd\ddlT \mathbf{x}_B\\
&\text{indeed }\left\{\begin{aligned}
  &\mathbf{x}_A = \left(\efassemp\ddl \impe\ddd\efassemp\ddlT \right)^{-1}\efassemp\ddl \efxi\ddc\\
  &\mathbf{x}_B = \left( \efassemd\ddl {\impe\ddd}^{-1}\efassemd\ddlT \right)^{+}\efassemd\ddl {\impe\ddd}^{-1}\efxi\ddc
  \end{aligned}\right.
\end{aligned}
\end{equation}
\subsection{Note on the existence of nonlinear Schur complements}
In previous sections we assumed the existence of nonlinear primal/dual/mixed Schur complement. For a well-posed global mechanical problem, local existence of such subdomain operator is most probable (see for instance \cite{ciarletLNL} for a list of nonlinear frameworks where solutions exist). Yet the formulation, the material and the shape of the subdomain may strongly limit the domain of existence and uniqueness of the Dirichlet/Neumann/Robin problems: typically the large displacements hypothesis may prevent us from using too large compressive loads because of the possibility of buckling, as well as damage may prevent us  from using too large tensile loads (see Section~\ref{section:assess} for an illustration of these difficulties).

The global existence of the operators can be proved in the case of coercive continuous monotone operators; see for instance \cite{showalter1977,showalter1997} for an analysis at the level of the variational formulation and \cite{Feistauer87} for the analysis of the finite element approximation. Mechanically this framework is associated with positive hardening behaviors and certain contact laws, in small strains \cite{Lad85,Lad99d}. In that context, Dirichlet, Neumann (assuming balanced load) and Robin problems are well-posed: the solution exists, is unique and varies continuously with respect to the load. Thanks to the continuity of the trace, the Schur operators can be defined. Moreover, the Schur complements inherit properties from the global problem (typically monotonicity, coercivity and continuity;  see for instance \cite{Hauer15} and associated bibliography), which makes the global condensed problems (\ref{eq:primal_condensed},\ref{eq:dual_condensed},\ref{eq:mixed_condensed}) well-posed.

\section{Solution strategy}\label{sec:solstrat}
At this point, assuming the existence of local ``Schur-type'' nonlinear operators on subdomains, we have obtained one global nonlinear interface problem with an additive structure since it can be written as the assembly of subdomain contributions:
\begin{equation}\label{eq:geneNL}
 \mathbf{L}(\efx):= \sum_s \mathbf{L}\s(\efx)=0
\end{equation}
We propose to apply a Newton-Raphson procedure to that system. Iteration $k$ consists in solving:
\begin{equation}\label{eq:geneNewt}
\begin{aligned}
\frac{\partial \mathbf{L}}{\partial \efx}(\efx_k)\overset{\circ}{\efx}_k   & = - \mathbf{L}(\efx_k) \quad \text{ and }\quad \efx_{k+1}=\efx_k+\overset{\circ}{\efx}_k 
\\ \text{i.e. }\quad
 \left(\sum_s \frac{\partial \mathbf{L}\s}{\partial \efx}(\efx_k)\right)\overset{\circ}{\efx}_k &=-\sum_s \mathbf{L}\s(\efx_k)\quad \text{ and }\quad \efx_{k+1}=\efx_k+\overset{\circ}{\efx}_k 
\end{aligned}
\end{equation}
We observe that the right-hand side, which corresponds to the evaluation of the residual, is obtained as the assembly of the result of independent solutions to nonlinear problems for each subdomain; whereas the left-hand side corresponds to the construction of the tangent operator as the assembly of subdomains' tangent operators.

The key point of this solution procedure is that the tangent condensed operators can actually be  computed: indeed the nonlinear condensed operators are defined as an assembly of the trace (which are linear operations) of the solution to subdomains' problems, the linearization of which leads to the definition of their tangent stiffness matrices. In other words, the tangent of the nonlinear condensed operator (primal, dual or mixed) is the condensation of tangent subdomain operators. \medskip

We now explain in more detail the computation of the tangent operator for the primal approach. Let $\schurp_{t}\s$ be the tangent operator of the nonlinear Schur complement $\schurp_{nl}\s$. It can be defined as the unique linear operator satisfying the following relation for any $\overset{\circ}{\efdep}_b\s$:
\begin{equation*}
\overset{\circ}{\efl}_b\s:= \schurp_{nl}\s(\efdepik\s+\overset{\circ}{\efdep}_b\s;\efforce_{ext}\s)-\schurp_{nl}\s(\efdepik\s;\efforce_{ext}\s)= \schurp_{t_k}\s \overset{\circ}{\efdep}_b\s+ o(\overset{\circ}{\efdep}_b\s)
\end{equation*}
The associated subdomain problems can be written as:
\begin{equation*}
\begin{aligned}\efforce_{int}\s(\efdep_k\s)+\efforce\s_{ext}+\eftrace\sT\eflamk\s=0 \text{ with } \eftrace\s \efdep_k\s = \efdepik\s \end{aligned}
\end{equation*}
The differentiation around $\efdep_k\s$ involves the tangent stiffness matrix $\efstiff_{t_k}\s$:
\begin{equation*}
\begin{aligned}-\efstiff_{t_k}\s\overset{\circ}{\efdep}\s + \eftrace\sT \overset{\circ}{\efl}_b\s + o\left(\overset{\circ}{\efdep}\s\right) = 0\text{ with } \eftrace\s \overset{\circ}{\efdep}\s = \overset{\circ}{\efdepi}\s 
\end{aligned}
\end{equation*}
If we assume the well-posedness of the tangent Dirichlet problem ($\|(\efstiff_{t_k ii}\s)^{-1}\|<\infty$), then we can condense the previous equation without amplifying the negligible terms and then we obtain that 
\begin{equation*}
\begin{aligned}
\overset{\circ}{\efl}_b\s & =\schurp_{t_k}\s \overset{\circ}{\efdep}_b\s+o\left(\overset{\circ}{\efdep}_b\s\right) \\
&=\left(\efstiff_{t_k bb}\s-\efstiff_{t_k bi}\s(\efstiff_{t_k ii}\s)^{-1}\efstiff_{t_k ib}\s\right)\overset{\circ}{\efdep}_b\s  + o\left(\overset{\circ}{\efdep}_b\s - \efstiff_{t_k bi}\s(\efstiff_{t_k ii}\s)^{-1}\overset{\circ}{\efdep}_i\s\right)
\end{aligned}
\end{equation*}
Hence the identification between the tangent Schur complement and the Schur complement of the tangent stiffness matrix.
\medskip

In the following, we derive the important steps of the different strategies.

\subsection{Primal formulation}
In the primal case, applying the Newton algorithm to equation \eqref{eq:primal_condensed} leads to:
\begin{equation}\label{eq:primal_newton}
\left(\efassemp\ddl   \schurp_{t_{k}}\ddd    \efassemp\ddlT\right){\ovefddepik} = - \efassemp\ddl\schurp_{nl}\ddc(\efassemp\ddlT{\ovefdepik};\efforce_{ext}\ddc)
\end{equation}
where $\schurp_{t_{k}}\s=\left(\frac{\partial\schurp_{nl}\s}{\partial \efdepi\s }(\efdepik\s;\efforce_{ext}\s)\right)$ is the subdomain  tangent primal Schur complement as given in \eqref{eq:primal_linear}.

First, let us focus on the right-hand side. It corresponds to the evaluation of a nonlinear Dirichlet problem for each subdomain with imposed interface displacements $\ovefdepik$:
\begin{equation}\label{eq:primform2}
\efdep_k\ddc \text{ such that }\left\{\begin{aligned}& \left(\efforce_{int}\ddc(\efdep_k\ddc)+\efforce\ddc_{ext}\right)_i=0\\
&\eftrace\ddd \efdep_k\ddc = \efassemp\ddlT \ovefdepik\end{aligned}\right.
,\quad \eflamk\ddc:=-\left(\efforce_{int}\ddc(\efdep_k\ddc)+\efforce\ddc_{ext}\right)_b
\end{equation}
From the parallel solution of these systems, we obtain the internal displacements $\efdep_k\ddc$ and the reactions $\eflamk\ddc$ whose lack of balance $\efassemp\ddl \eflamk\ddc=\efassemp\ddl\schurp_{nl}\ddc(\efassemp\ddlT{\ovefdepik};\efforce_{ext}\ddc)$ is the right-hand side of the tangent system.

The left-hand side of \eqref{eq:primal_newton} is the assembly of the subdomain's tangent primal Schur complements, that is to say a classical primal domain decomposition formulation of the tangent problem. 

Thus, after an arbitrary initialization, the straightforward solution consists in repeating the following steps: (i) solving independent  nonlinear Dirichlet systems for each subdomain, (ii) computing the lack of balance at the interface, and (iii) solving the global tangent interface problem by the application of the BDD algorithm \cite{Man93,Let94}.

\subsection{Dual formulation}
In the dual case, applying Newton algorithm to equation \eqref{eq:dual_condensed} leads to:
\begin{equation}\label{eq:dual_newton}
\left\{\begin{aligned}
&\left(\efassemd\ddl   \schurd_{t_k}\ddd    \efassemd\ddlT\right){\uvefdlamk}+\efassemd\ddl\efker_b\ddd\overset{\circ}{\efalp}_k \ddc =  - \efassemd\ddl\left(\schurd_{nl}\ddc(\efassemd\ddlT{\uveflamk};\efforce_{ext}\ddc)+\efker_b\ddd\efalp_k\ddc\right)\\
&\efker_b\dddT\efassemd\ddlT {\uvefdlamk}=0
\end{aligned}\right.
\end{equation}
where $\schurd_{t_k}\s=\left(\frac{\partial\schurd_{nl}\s}{\partial \eflam\s }(\eflamk\s;\efforce_{ext}\s)\right)$ is the subdomains' tangent dual Schur complement in configuration $k$ as given in \eqref{eq:dual_linear}.

As usually done in the FETI method, we introduce a projector $\proj$ on $\ker(\efker_b\dddT\efassemd\ddlT )$, so that we can seek $\proj \uvefdlamk$ instead of $\uvefdlamk$, solution to:
\begin{equation}
\proj^T\left(\efassemd\ddl   \schurd_{t_k}\ddd    \efassemd\ddlT\right)\proj{\uvefdlamk} =  - \proj^T\efassemd\ddl\schurd_{nl}\ddc(\efassemd\ddlT{\uveflamk};\efforce_{ext}\ddc)\\
\end{equation}
and the contribution of the rigid body motions is sought after $\uveflam$ is determined.

The right-hand side of \eqref{eq:dual_newton} corresponds to the evaluation of nonlinear Neumann problems for each subdomain with imposed interface reaction $\uveflamk$:
\begin{equation}\label{eq:dualform2}
 \efforce_{int}\ddc(\efdep_k\ddc)+\efforce\ddc_{ext}+\eftrace\dddT\efassemd\ddlT{\uveflamk}=0
\end{equation}
From the parallel solution of these systems, we obtain the displacements $\efdep_k\ddc$ whose interface gap  $\proj^T\efassemd\ddl {\eftrace\ddd} \efdep_k\ddc=\proj^T\efassemd\ddl\schurd_{nl}\ddc(\efassemd\ddlT{\uveflamk};\efforce_{ext}\ddc)$ is the right-hand side of the tangent system.

The tangent matrix is the assembly of local tangent dual Schur complements computed at the value $\efdep\ddc_k$ obtained when computing the right-hand side. Then the first term of \eqref{eq:dual_newton} is the assembly of the subdomain's tangent dual Schur complements, that is to say a classical dual domain decomposition formulation of the tangent problem. 

In order for system \eqref{eq:dualform2} to be well posed, the rigid-body admissibility condition \eqref{eq:admiss} has to be satisfied. This is realized by initializing $\uveflam$ by an admissible reaction ${\uveflamz}$ which verifies:
\begin{equation}
\efker\dddT\left(\efforce_{ext}\ddc+\eftrace\dddT \efassemd\ddlT {\uveflamz}\right)=0\end{equation}
 Indeed, since increments ${\uvefdlamk}$ satisfy a zero-admissibility condition \eqref{eq:dual_newton}, and because the admissibility is a linear condition thanks to the small perturbation hypothesis, any coming $\uveflamk$ will satisfy the admissibility condition. Finding an initial value $\uveflamz$ corresponds to the initialization of the FETI method.

Thus, the straightforward solution consists in solving a FETI coarse grid problem to initialize $\uveflamz$, then repeating the following steps: (i) solving independent   nonlinear Neumann systems for each subdomain, (ii) computing the displacement gap at the interface, and (iii) solving the global tangent interface problem by the application of FETI algorithm \cite{Far94}. Note that the initialization can be improved by a full solution of the initial FETI system.

\subsection{Mixed formulation}
In the mixed case applying Newton algorithm to equation \eqref{eq:mixed_condensed} leads to:
\begin{multline}\label{eq:mixed_newton}
  \left(\efassemp\ddlT\left( \efassemp\ddl \impe\ddd\efassemp\ddlT \right)^{-1} \efassemp\ddl   - \schurm\ddd_{t_k}\right) \efdmuk\ddc \\= \schurm\ddc_{nl}(\efmuk\ddc;\efforce_{ext}\ddc,\impe\ddd)-\efassemp\ddlT\left( \efassemp\ddl \impe\ddd\efassemp\ddlT \right)^{-1} \efassemp\ddl\efmuk\ddc
\end{multline}
where $\schurm_{t_k}\s= \left(\frac{\partial\schurm_{nl}\s}{\partial \efmu\s }(\efmuk\s;\efforce_{ext}\s,\impe\s)\right)$ is the subdomains' tangent mixed Schur complement in configuration $k$ as given in \eqref{eq:mixed_linear}.

The first term of the right-hand side of \eqref{eq:mixed_newton} corresponds to the evaluation of nonlinear Robin problems for each subdomain with interface impedance $\impe\ddd$ and imposed interface reaction $\efmuk\ddc$:
\begin{equation}\label{eq:mixedform2}
\efforce_{int}\ddc(\efdep_k\ddc)-\eftrace\dddT\impe\ddd\eftrace\ddd\efdep_k\ddc+\eftrace\dddT\efmuk\ddc+\efforce\ddc_{ext}=0
\end{equation}
from the parallel solution of these systems, we obtain the displacement $\efdep_k\ddc$ from which we deduce the mixed residual $\efrhs_{m_k}\ddc=\efdepik\ddc-\efassemp\ddlT\left( \efassemp\ddl \impe\ddd\efassemp\ddlT \right)^{-1} \efassemp\ddl\efmuk\ddc$ which is the right-hand side of the tangent system.

The tangent operator is obtained by mixed condensation of the tangent stiffness matrix computed at the value $\efdep\ddc_k$ obtained when evaluating the right-hand side. Then the first term of \eqref{eq:mixed_newton} is a classical mixed domain decomposition formulation of the tangent problem. 

Thus after an arbitrary initialization, the straightforward solution consists in repeating the following steps: (i) solving independent  nonlinear Robin systems for each subdomain, (ii) computing the mixed residual, and (iii) solving the global tangent interface problem by the application of FETI-2LM algorithm \cite{FXFETI2LM09}.

\subsection{Swapping the linear solvers}
In any case, the solution procedure alternates independent nonlinear solutions for each subdomain with different boundary conditions, and global tangent solutions similar to a domain decomposition method applied to linear problems. As we will show later, the choice of the boundary condition strongly impacts the nonlinear subdomains' solving. As the tangent operator is always associated with the tangent stiffness matrix of subdomains, it is possible to choose any formulation for the linear problem.

Indeed, starting from (\ref{eq:primal_linear},\ref{eq:dual_linear},\ref{eq:mixed_linear}) we have the following classical relationships~\cite{GOSSELET.2007.1}:
\begin{equation}
  \begin{aligned}
    \schurd_t\ddd&={\schurp_t\ddd}^+            &\qquad \schurm_t\ddd&=\left(\schurp_t\ddd+\impe\ddd\right)^{-1} \\
    \efrhs_d\ddc &=  \schurd_t\ddd \efrhs_p\ddc &\qquad \efrhs_m\ddc &=   \schurm_t\ddd \efrhs_p\ddc
  \end{aligned}
\end{equation}
This implies that even if a formulation results in a specific residual after the subdomains' nonlinear solutions, any linear solver can be employed. Typically, the mixed formulation has interesting capabilities for nonlinear problems, but it is not very classical for linear system for which primal (BDD) or dual (FETI) formulations are more standard and embed powerful preconditioners. After computing $\schurm\ddc_{nl}(\efmuk\ddc;\efforce_{ext}\ddc,\impe\ddd)$, one can compute $\efrhs_{m_k}\ddc = \schurm\ddc_{nl}(\efmuk\ddc;\efforce_{ext}\ddc,\impe\ddd)-\efassemp\ddlT\left( \efassemp\ddl \impe\ddd\efassemp\ddlT \right)^{-1} \efassemp\ddl\efmuk\ddc$, deduce the equivalent primal right-hand side $\efrhs_{p_k}\ddc=\left(\schurp_{t_k}\ddd+\impe\ddd\right)\efrhs_{m_k}\ddc$ and use the associated interface lack of balance $\efassemp\ddl\efrhs_{p_k}\ddc$ as input to the tangent BDD solver. This property can directly be observed by introducing \eqref{eq:mixedortho3} in the  mixed tangent system \eqref{eq:mixed_newton}, since after premultiplying by $(\impe\ddd+\schurp_{t_k}\ddd)$ we obtain:
\begin{equation}\label{eq:syseqmixt}
\begin{aligned}
  \schurp_{t_k}\ddd \efassemp\ddlT\ovefddepiik - \efassemd\ddlT\uvefdlammk &=  (\impe\ddd+\schurp_{t_k}\ddd)\left(\schurm\ddc_{nl}(\efmuk\ddc;\efforce_{ext}\ddc,\impe\ddd)-\efassemp\ddlT \ovefdepiik\right)\\
  &=(\impe\ddd+\schurp_{t_k}\ddd) \efrhs_{m_k}\ddc = \efrhs_{p_k}\ddc
\end{aligned}
\end{equation}
which leads to the system set in terms of $\ovefddepiik$:
\begin{equation}
\begin{aligned}
\left( \efassemp\ddl \schurp_{t_k}\ddd \efassemp\ddlT \right) \ovefddepiik  &=  \efassemp\ddl\efrhs_{p_k}\ddc 
\end{aligned}
\end{equation}
The following quantities can be deduced (even from an inexact solution):
 \begin{equation}\label{eq:fromdv}
\begin{aligned}
\efdmuk\ddc &=  \schurp_{t_k}\ddd \efassemp\ddlT \ovefddepiik-\efrhs_{p_k}\ddc \\
\efddepk\ddc&=(\efstiff_{t_k}\ddd+\eftrace\dddT\impe\ddd\eftrace\ddd)^{-1}\eftrace\dddT\left(\efrhs_{p_k}\ddc+\efdmuk\ddc\right)\\
\efddepik\ddc&=\eftrace\ddd\efddepk\ddc
\\ \efdlamk\ddc & = \schurp_{t_k}\ddd \efddepik\ddc - \efrhs_{p_k}= \efdmuk\ddc - \impe\ddd\efddepik\ddc
\end{aligned}
\end{equation}
Note that $\uvefdlammk$ plays no role in the algorithm but it could also be post-processed:
\begin{equation}
\begin{aligned}
\uvefdlammk &=  {\efassemd\ddlT}^+\left( \efdmuk\ddc -\impe\ddd \efassemp\ddlT \ovefddepiik\right)
\end{aligned}
\end{equation}
where ${\efassemd\ddlT}^+$  is a scaled assembling operator commonly employed in the preconditioning step of the FETI method \cite{Rix99bis,Kla01}; because of the space splitting properties \eqref{eq:mixedortho2}, $\uvefdlammk$ does not depend on the choice of the scaling. Of course, if a dual formulation had been preferred then $\uvefdlammk$ would have been computed by a FETI-like system (left-multiply \eqref{eq:syseqmixt} by $\efassemd\ddl\schurd_t\ddd$).

\subsection{Typical algorithm}

Algorithm~\ref{alg:robin-bdd} sums up the main steps of the method with the mixed nonlinear local problems and primal global solver. For simplicity reasons, only one increment was considered.

As can be seen in this algorithm, besides classical global convergence criterion $\varepsilon_{NG}$ (as in a standard Newton approach), two precision thresholds are used: the local nonlinear thresholds $\varepsilon_{NL}\ddc$ (associated with the Newton processes carried out independently on subdomains) and the global linear threshold of the domain decomposition (Krylov)  solver $\varepsilon_{K}$ (here BDD).
The convergence criteria are discussed in more  detail in the following section.

 Other parameters can be the initializations of the various iterative solvers and, for the mixed approach, the choice of the impedance matrices $\impe\ddd$. 

\begin{algorithm2e}\caption{Mixed nonlinear approach with BDD tangent solver}\label{alg:robin-bdd}
\dontprintsemicolon
\KwSty{Define:}\;
$\efres_{nl}^{m\diamondvert}(\efdep\ddc,\efmu\ddc) =\efforce_{int}\ddc(\efdep\ddc)-\eftrace\dddT\impe\ddd\eftrace\ddd\efdep\ddc+\eftrace\dddT{\efmu}\ddc+\efforce\ddc_{ext}$\;
\BlankLine
\KwSty{Initialization:}\;
$(\efdep_0\ddc,{\eflam}_0\ddc)$ such that $\efassemd\ddl\eftrace\ddd\efdep_0\ddc=0$ and $\efassemp\ddl{\eflam}_0\ddc=0$\;
\KwSty{Set} $k=0$\;
\KwSty{Define} ${\efmu}_k\ddc={\eflam}_{k}\ddc + \impe\ddd \eftrace\ddd{\efdep}_{k}\ddc$\;
\While{$\|\efres_{nl}^{m\diamondvert}(\efdep_k\ddc,{\efmu}_k\ddc)\|+ \|\efassemd\ddl\eftrace\ddd\efdep\ddc\|_{B}>\varepsilon_{NG}$ }{%
  \KwSty{Local nonlinear step}:  \;
   \KwSty{Set} ${\efdep}_{k,0}\ddc={\efdep}_{k}\ddc$ and  $j=0$\; 
  \While{$\|\efres_{nl}^{m\diamondvert}(\efdep_{k,j}\ddc,{\efmu}_k\ddc)\|>\varepsilon\ddc_{NL}$}{
    ${\efdep}_{k,j+1}\ddc={\efdep}_{k,j}\ddc-\left({\stiff_t}_{k,j}\ddd+\eftrace\dddT\impe\ddd\eftrace\ddd\right)^{-1}\efres_{nl}^{m\diamondvert}(\efdep_{k,j}\ddc,{\efmu}_k\ddc)$ \;
    \KwSty{Set} $j=j+1$
  }
  \KwSty{Linear right-hand side}:\;
  ${\efrhs_m}_k\ddc = \efassemp\ddlT\left( \efassemp\ddl \impe\ddd\efassemp\ddlT \right)^{-1}    \efassemp\ddl  {\efmu}_k\ddc -{\eftrace\ddd}{\efdep}_{k,j}\ddc $\;
  ${\efrhs_p}_k\ddc = ({\schurp_t}_{k,j}\ddd +\impe\ddd){\efrhs_m}_k\ddc $\;
  \KwSty{Global linear step}:\;
  \KwSty{Set} $\ovefddepiik^{0}=0$ and  $i=0$\; 
  \While{$\|{\efrhs_p}\ddc_k-\left(\efassemp\ddl\; {\schurp_t}_{k,j}\ddd\efassemp\ddlT\right)\ovefddepiik^i\|>\varepsilon_{K}$}{
    Make BDD iterations (index $i$)
  }
  \KwSty{Set} ${\efdep}_{k+1}\ddc = \efdep_k\ddc + {\overset{\circ}{\efdep}_k^{i\diamondvert}}$ and ${\eflam}_{k+1}\ddc = \eflamk\ddc + \efdlamk^{i\diamondvert}$ using \eqref{eq:fromdv}\;
  \KwSty{Set} $k=k+1$\;
}
\end{algorithm2e}

\section{Error analysis}\label{sec:error}
The aim of this section is to further analyze the effect of the different thresholds and the way they are connected. When possible, we will make reference to the global residual $\efres$ which is defined as:
\begin{equation}\label{eq:residglocal}
\efforce_{int}(\efdep)+\efforce_{ext} = \efres 
\end{equation}
The objective of the solver is that $\|\efres\|<\varepsilon_{NG}\|\efforce_{ext}\|$ (global nonlinear criterion). 

We first show how the convergence criteria for the Local Newton solvers $\varepsilon\ddc_{NL}$ and for the  Krylov interface solver $\varepsilon_K$ are related for the different approaches. Then we interpret the whole  process as an inexact Newton solver \cite{dembo1982inexact} and we consider adapting at each step the criteria $\varepsilon_{K}$ and $\varepsilon_{NL}\ddc$  to the current error status $\|\efres\|$.  $\varepsilon_{K}$ and $\varepsilon_{NL}\ddc$ are interpreted as forcing terms which need to be close to the objective $\varepsilon_{NG}$ only when global convergence is almost reached whereas they should be relaxed at the beginning of the process in order to avoid oversolving (see Section~\ref{section:sync}).

\subsection{Primal approach}
\subsubsection{Local Newton solvers}
The inner (local) Newton loops are associated with Dirichlet problems. Assuming these conditions are exactly taken into account, the convergence is controlled by the internal node residue:
\begin{equation}\label{eq:residprimallocal}
{\efforce\ddc_{int}}_i(\efdep\ddc)+{\efforce\ddc_{ext}}_i = \efres\ddc_i\text{ with }\|\efres_i\ddc\|<\varepsilon_{NL}\ddc\|{\efforce\ddc_{ext}}_i\|
\end{equation}
From the $\efdep\ddc$ computed above, the reactions are defined by the following relation:
\begin{equation}\label{eq:deflam}
\eflam\ddc:=-\left({\efforce_{int}}_b\ddc(\efdep\ddc)+{\efforce\ddc_{ext}}_b\right)
\end{equation}

\subsubsection{Link to the global error}
In the primal case, the displacement search space is the same as in the non-substructured case (displacements are continuous across the interface). This simplifies the error analysis. Indeed, we directly have:
\begin{equation}\label{eq:primal_split}
\begin{aligned}
\|\efforce_{int}(\efdep)+\efforce_{ext}\| &\leqslant \|\efforce_{int_i}\ddc(\efdep\ddc)+\efforce_{ext_i}\ddc\| + \|\efassemp\ddl\left(\efforce_{int_b}\ddc(\efdep\ddc)+\efforce_{ext_b}\ddc\right)\|\\
&\leqslant \|\varepsilon_{NL}\ddc\| \|\efforce_{ext_i}\ddc\|+ \|\efassemp\ddl\eflam\ddc\|
\end{aligned}
\end{equation}
This means that global convergence occurs when both the local Newton criteria and the resulting lack of balance at the interface (the initial residual of the tangent interface solver) are small.

\subsubsection{Tangent interface solver}
From the point of view of the nonlinear Schur complement, the computation of the reactions suffers from the error made on internal displacements, if $\schurp_{nl}\ddc(\efassemp\ddlT{\ovefdepik};\efforce_{ext}\ddc)$ denotes the exact reaction (corresponding to $\efres\ddc_i=0$), we write:
\begin{equation}\eflam\ddc=\schurp_{nl}\ddc(\efassemp\ddlT{\ovefdepik};\efforce_{ext}\ddc)+\efres_b\ddc
\end{equation}
 In the linear case, the internal residual would propagate on the boundary according to the relation $\efres_b\ddc=-\stiff_{bi}\ddd{\stiff_{ii}\ddd}^{-1}\efres_i\ddc$, so that we can assume that the error on the reaction is of the same order of the inner solver error $\|\efres_b\ddc\|\simeq\|\efres_i\ddc\|$.

We then can study the error made on the tangent solution. The computed vector $\ovefddepik$ satisfies the equation:
\begin{equation}
\left(\efassemp\ddl   \schurp_{t_{k}}\ddd    \efassemp\ddlT\right){\ovefddepik} = - \efassemp\ddl\schurp_{nl}\ddc(\efassemp\ddlT{\ovefdepik};\efforce_{ext}\ddc) - \efassemp\ddl \efres_b\ddc + \efres_{\!A}^{L}
\end{equation}
where $\|\efres_{\!A}^{L}\|<\varepsilon_K \|\efres_{\!A_0}^{L}\|$ controls the error of the linear (Krylov) solver. Since $\|\efassemp\ddl\|=O(1)$, we see that the error of the local Newton solver and the global interface iterative solver should be of the same order of magnitude.

\subsection{Dual approach}
\subsubsection{Local Newton solvers}
The inner (local) Newton loops are associated with Neumann problems, and the convergence is controlled by the following residue:
\begin{equation}
\efforce_{int}\ddc(\efdep\ddc)+\efforce\ddc_{ext}+\eftrace\dddT\efassemd\ddlT\uveflam=\efres\ddc \text{ with } \|\efres\ddc\|<\varepsilon_{NL}\ddc\|\efforce\ddc_{ext}\|
\end{equation}
The output from these computations is the interface displacement. Unfortunately there is no direct control on the associated error $\boldsymbol{\delta}_b\ddc$. A good estimation is provided by the last correction brought to $\efdep_b\ddc$ by the inner Newton loop.
 
\subsubsection{Link to the global error}
The difficulty with the dual approach is that the (broken) search space is larger than the original space (where displacements are continuous). In the linear case, a costless processing enables us to obtain a continuous displacement \cite{PARRETFREAUD.2010.1.1}.  Such a strategy is not realistic for nonlinear problems.

One possibility is to define an error extended to the broken space ($e_{bs}$):
\begin{equation}\label{eq:residglodual}
e_{bs}^2=\|\efforce_{int}\ddc(\efdep\ddc)+\efforce_{ext}\ddc+\eftrace\dddT\efassemd\ddlT{\uveflam}\|^2 + \|\efassemd\ddl\eftrace\ddd\efdep\ddc\|_{B}^2 
\end{equation}
where $\| \|_{B}$ is a well chosen norm for the interface discontinuity. We recover an expression where the inner criterion is cumulated with the initial residual of the tangent interface solver.

\subsubsection{Tangent interface solver}
The tangent interface solution can be written as:
\begin{equation}
\proj^T\left(\efassemd\ddl   \schurd_{t_{k}}\ddd    \efassemd\ddlT\right)\proj{\uvefdlamk} = - \proj^T\left(\efassemd\ddl\schurd_{nl}\ddc(\efassemd\ddlT{\uveflamk};\efforce_{ext}\ddc) + \efassemd\ddl\boldsymbol{\delta}_b\ddc\right) + \efres_{\!B}^{L}
\end{equation}
where $\|\efres_{\!B}^{L}\|<\varepsilon_{K}\|\efres_{\!B_0}^{L}\|$ controls the error of the linear solver. Since $\|\efassemd\ddl\|=O(1)$, we see that the forward error of the local Newton solver and the backward error of the global interface iterative solver should be of the same order of magnitude. 

\subsection{Mixed approach}
\subsubsection{Local Newton solvers}
The inner (local) Newton loops are associated with Robin problems. The convergence is controlled by the following residue:
\begin{equation}
\efforce_{int}\ddc(\efdep\ddc)-\eftrace\dddT\impe\ddd\eftrace\ddd\efdep\ddc+\eftrace\dddT{\efmu}\ddc+\efforce\ddc_{ext}=\efres^{\diamondvert} \text{ with } \|\efres\ddc\|<\varepsilon_{NL}\ddc\|\efforce\ddc_{ext}\|
\end{equation}
As in the dual approach, the output from these computations is the interface displacement, of which the error $\boldsymbol{\delta}_b\ddc$ can be estimated by the last correction brought to $\efdep_b\ddc$ by the inner Newton loop.
\subsubsection{Link to the global error}
As in the dual case, the search space is a broken space. This implies the use of an extended norm to evaluate the global error: not only the displacements are not continuous at the interface but neither are the reactions balanced. With reaction $\eflam\ddc$ defined by \eqref{eq:deflam}, we see that the trace of the local residue can be written as:
\begin{equation}
\efres\ddc_b=\impe\ddd\left(\efassemp\ddlT \mathbf{v}_A-\eftrace\ddd\efdep\ddc\right)+ \left( \efassemd\ddlT \boldsymbol{\gamma}_B  - \eflam\ddc   \right)
\end{equation}
This relation links the two interface errors, which means that the dual criterion \eqref{eq:residglodual} is sufficient to monitor the convergence: if both $\|\efres\ddc_b\|$ and $\|\efassemd\ddl\eftrace\ddd\efdep\ddc\|_B$ are small, then so must $\|\efassemp\ddl\eflam\ddc\|$ be.
\subsubsection{Tangent interface solver}
The tangent interface solution can be written as (in primal form):
\begin{equation}
\left( \efassemp\ddl \schurp_{t_k}\ddd \efassemp\ddlT \right) \ovefddepiik  =  \efassemp\ddl(\impe\ddd+\schurp_{t_k}\ddd)\left(\efdepik\ddc+\boldsymbol{\delta}_b\ddc-\efassemp\ddlT \ovefdepiik\right)+\efres_{\!A}^L
\end{equation}
where $\|\efres_{\!A}^{L}\|<\varepsilon_{K}\|\efres_{\!A_0}^{L}\|$ controls the error of the linear solver. From the definition of $\boldsymbol{\delta}_b\ddc$, we can expect  $\|\efassemp\ddl(\impe\ddd+\schurp_{t_k}\ddd)\boldsymbol{\delta}_b\ddc\| \simeq \|\efassemp\ddl \eftrace\ddd\efres\ddc\|$ which means that the convergence criterion of the Krylov solver must be of the order of magnitude of the local Newton solvers' criteria.

\subsection{Tuning of criteria} \label{section:sync}
The previous subsections established connections between the convergence criteria of the local Newton solvers and of the global interface Krylov solver. We  can now analyze the algorithm in the framework of inexact Newton methods \cite{dembo1982inexact}. 
From previous analysis, and using the notations of \eqref{eq:geneNewt}, we solve \eqref{eq:geneNL} by the sequence:
\begin{equation}\label{eq:inexNewt}
\begin{aligned}
\frac{\partial \mathbf{L}}{\partial \efx}(\efx_k)\overset{\circ}{\efx}_k &=-\mathbf{L}(\efx_k)+ \efres_k,\text{ with }\|\efres_k\|\leqslant \varepsilon_k \|\mathbf{L}(\efx_k) \| \\
  \efx_{k+1}&=\efx_k+\overset{\circ}{\efx}_k 
\end{aligned}
\end{equation}
$\varepsilon_k$ is called the forcing term, and its value can be updated at each Newton iteration in order to optimize the convergence rate. More precisely,  \cite{eisenstat1996choosing} shows that:
\begin{enumerate}
\item If $0 \leq \varepsilon_k \leq \varepsilon_{max} < 1 \, \, \, \forall k$, then the series converges linearly, with an asymptotic rate $a \leqslant \varepsilon_{max}$.
\item If $\lim\limits_{k \rightarrow +\infty} \varepsilon_{k} = 0$, convergence is superlinear.
\item If $\varepsilon_k = O( \Vert\mathbf{L}(\efx_k) \Vert )$, convergence is quadratic.
\end{enumerate}

These convergence properties lead to more or less simple expressions for the \textit{forcing terms} (see \cite{eisenstat1996choosing}):
\begin{align}
\varepsilon_k & =\dfrac{1}{2^{k+1}} \label{expr:forcterm1} \\
\varepsilon_k & =\gamma \left( \frac{\Vert \mathbf{L}(\efx_k)\Vert}{\Vert  \mathbf{L}(\efx_{k-1})\Vert} \right)^{\alpha}, \, \, \gamma \, \in \, \left[ 0,1 \right[, \, \, \alpha \, \in \, \left] 0,2 \right] \label{expr:forcterm2} \\
\varepsilon_k & =\frac{\vert \Vert  \mathbf{L}(\efx_k)\Vert - \Vert \mathbf{L}(\efx_{k-1})+\mathbf{L}'(\efx_{k-1})\overset{\circ}{\efx}_{k-1} \Vert \vert}{\Vert  \mathbf{L}(\efx_{k-1})\Vert}\label{expr:forcterm3}
\end{align}

Expression~\eqref{expr:forcterm1} is easy to implement but does not relate the linear criterion with $\Vert \mathbf{L}(\efx_k) \Vert$ at each iteration: the convergence of the series $( \efx_k )_k$ will theoretically not be better than superlinear. Expressions~\eqref{expr:forcterm2} and \eqref{expr:forcterm3} are more complicated to implement, but can a priori allow faster convergence rates. For example taking $\alpha=2$ in expression~\eqref{expr:forcterm2} should lead to quadratic convergence.

\section{Assessments} \label{section:assess}
The aim of these assessments is to test the various configurations of the approach to an academic problem with a realistic nonlinear behavior.
We are aware that the methods we propose are sensitive to many parameters, among others the type of nonlinearities (with potential instabilities), the partitioning of the domain, the localization of the linearities (whether the nonlinearity spreads amongst subdomains or not), the variants of the solvers, the loading increments, the convergence criteria. 

The following study is thus far from pretending to be exhaustive, it simply aims at giving major trends of the methods and at proving that they  can be computationally interesting. At the end of the section, we  discuss  the applicability of the methods to other cases and the expected performance.

\subsection{Perfect plasticity problem}
\begin{figure}[ht]
\centering
\includegraphics[width=.6\textwidth]{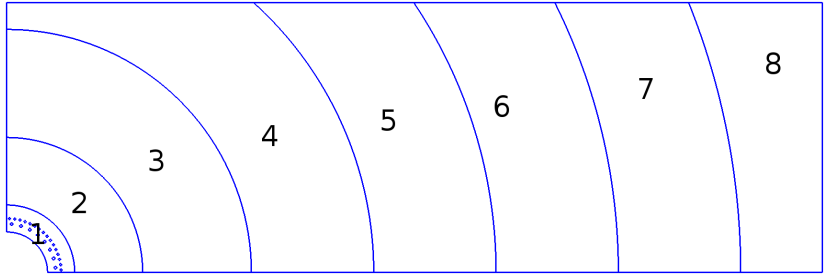}
\caption{Localized plasticity problem / numbering of subdomains}\label{fig:plast}
\end{figure}

We consider a rectangular plate in traction with a central hole (length $L_x=95$~mm, height $L_y=31.5$~mm, radius $r=5$~mm) as depicted in Figure~\ref{fig:plast} (after taking into account the symmetries).
The material is a perfectly plastic steel (Young modulus $E=210$~GPa, Poisson's ratio $\nu=0.3$, elasticity limit $\sigma_0=420$~MPa), and we assume the plane stress conditions. The mesh is constituted by triangle elements with three nodes, with characteristic length $h=0.01 \, L_y$. 
The eight subdomains are delimited by concentric circles. They are numbered starting from the center. 

Small perforations were added near the central hole in order to trigger a large level of localized plasticity. 
Traction loading is incrementally applied on the perforated plate by imposing a displacement $u_D$ on its right side. Increments are taken as follows, with $u_{e} = 25$~$\mu$m:
\begin{equation*}
u_D=\left[ 1 \quad 2 \quad 3 \quad 4 \quad 5 \quad 5.75 \quad 6.5 \right] u_{e}
\end{equation*}
The first step, corresponding to an increment of loading $u_D=u_{e}$, remains elastic. Steps $u_D=2 u_{e}$ to $u_D=5 u_{e}$ activate plasticity in only the central subdomain, reaching a maximum of $p=0.051$ of cumulated plasticity. Plasticity then spreads in the second and third subdomains until $\max(p)=0.15$.

This sequence of loadings was identified using a classical (outer Newton) approach. It could be reused as is with primal and mixed approaches. Solving a local Neumann problem with perfectly plastic subdomains is known to be a tougher problem than a Dirichlet or a Robin problem. Indeed, it was necessary to adopt a finer loading sequence for the dual approach; the following steps needed to be inserted:
\begin{equation*}
u_D = [ 4.5 \quad 5.25 \quad 5.5 \quad 6 \quad 6.125 \quad 6.25 \quad 6.375 ] u_e
\end{equation*}
It was verified that, for the primal and mixed approaches, this refined sequence leads to results of comparable quality as the coarse sequence.

In the following, we compare the classical approach (outer Newton, inner Schur-Krylov), the primal and dual nonlinear approaches, and the mixed nonlinear approach with parameter $\impe$ equal to the stiffness of the neighbors, or $\impe$ equal to the optimal choice (Schur complement of the remainder of the structure). Of course the latter choice is not computationally realistic, and efficient approximations will be the subject of further studies.

In order to assess the performance, we measure:
\begin{itemize}
\item The number of outer Newton loops, since nonlinear approaches shall enable us to follow a faster path to convergence, 
\item The maximum number of computations of the tangent amongst subdomains. For the classical approach, it is simply equal to the number of outer iterations. The inner Newton loops makes that number larger in non-linear approaches.
\item The number of cumulated Krylov (conjugate gradient) iterations which is proportional to the number of exchanges between subdomains. 
\end{itemize}
The aim of the nonlinear approaches is to lower the number of outer Newton and Krylov iterations, at the price of some extra local Newton iterations.   Of course the number of computations of the tangent varies from one subdomain to another, depending on the spread and intensity of the nonlinearity, which can cause poor load balancing. This is a real issue which will be the subject of future studies. Nevertheless we believe that the reduction of the number of outer iterations and of the exchanges (Krylov iterations) already makes the proposed methods interesting.

\subsubsection{Classical technique}
We briefly give the performance for the classical Newton method with a domain decomposition solver for the tangent systems. Note that BDD is used but FETI gives a quite similar number of Krylov iterations. The outer Newton's convergence criterion is set to $\varepsilon_{NG} = 10^{-5}$. For the results given in Table \ref{res:classical}, the inner Krylov solver criteria for the condensed system is fixed to $\varepsilon_{K} = 10^{-6}$.

\begin{table}[ht]
\begin{center}
\renewcommand{\arraystretch}{1.5}
\begin{tabular}{|m{4cm}|c|c|c|c|c|}
\hline
Increment of loading & {$u_e$} &{2$u_e$} & {5$u_e$} & {5.75$u_e$} & {6.5$u_e$} \\
\hline
Spread of nonlinearity & {Elastic} &  \multicolumn{2}{c|}{1SD plastifies} & \multicolumn{2}{c|}{Several SD plastify}\\
\hline
Global iterations           & 1 & 2 & 13 & 17 & 23\\\hline
Cumulated Krylov iterations & 19 & 38 & 255 & 336 & 460 \\\hline
\end{tabular}
\end{center}
\caption{Performance of the classical approach with BDD solver for the tangent systems, with fixed criteria.}
\label{res:classical}
\end{table}

\subsubsection{Nonlinear localization with fixed criteria}

Before tuning the different thresholds, as developed in section~\ref{section:sync}, a comparison is made between the classical outer Newton and the nonlinear approaches with fixed stopping criteria. In order to achieve a meaningful comparison, thresholds are fixed to typical values:
\begin{center}
$\varepsilon_{NG} = 10^{-5}$\\
$\varepsilon_{NL} = \varepsilon_{K} = 10^{-6}$
\end{center}

Relative results of primal, dual and mixed approaches normalized by the classical approach are given in Table \ref{res:global_fix} for chosen loading increments. 

Table \ref{res:global_fix} clearly indicates a decrease in the number of cumulated global iterations when using the primal and mixed approaches: 20\% less iterations for the primal and $\stiff_{bb}$-mixed approaches, 30\% less iterations with the optimal Robin condition. For the dual method, the number of global iterations increases because of the steps inserted additionally to insure convergence.

\begin{table}[ht]
\begin{center}
\renewcommand{\arraystretch}{1.5}
\begin{tabular}{|c|m{3.8cm}|c|c|c|c|c|}
\hline
&Increment of loading & {$u_e$} &{2$u_e$} & {5$u_e$} & {5.75$u_e$} & {6.5$u_e$} \\
\hline
&Spread of nonlinearity & {Elastic} &  \multicolumn{2}{c|}{SD 1 } & SD 1-2 & SD 1-3\\
\hline
\multirow{4}{*}{\rotatebox[origin=c]{90}{Global Newton}}&Primal/Classic & 1 & 1 & 0.77 & 0.82 & 0.83 \\
&Dual/Classic & 1 & 1 & 1.46 & 2.12 & 3.35 \\
&Mixed/Classic, $\impe=\stiff_{bb}$ & 1 & 1 & 0.77 & 0.82 & 0.78 \\
&Mixed/Classic, $\impe$ opti & 1 & 1 & 0.62 & 0.65 & 0.70 \\
\hline
\multirow{4}{*}{\rotatebox[origin=c]{90}{Krylov}}
&Primal/Classic & 1 & 1 & 0.76 & 0.82 & 0.82 \\
&Dual/Classic & 0.95 & 0.95 & 1.36 & 1.96 & 3.08\\
&Mixed/Classic, $\impe=\stiff_{bb}$ & 1 & 1 & 0.76 & 0.82 & 0.78 \\
&Mixed/Classic, $\impe$ opti & 1 & 1 & 0.61 & 0.64 & 0.69 \\\hline
\multirow{4}{*}{\rotatebox[origin=c]{90}{Local Newton}}
&Primal/Classic & 1 & 2 & 2.62 & 2.88 & 3.04 \\
&Dual/Classic & 1 & 2 & 5.69 & 9.24 & 16.74 \\
&Mixed/Classic, $\impe=\stiff_{bb}$ & 1 & 2 & 2.69 & 2.94 & 3 \\
&Mixed/Classic, $\impe$ opti & 1 & 2 & 1.85 & 2.12 & 2.43 \\
\hline
\end{tabular}
\end{center}
\caption{Ratios of cumulated iterations between nonlinearly localized methods and classic method, with fixed stopping criteria}
\label{res:global_fix}
\end{table}

The gain in Krylov iterations is quite similar to the gain in global Newton iterations, because each solution of a linear system requires an almost constant number of Krylov iterations (note that before requiring extra increments, the dual approach needs one Krylov iteration less than the primal approach, hence the $0.95\simeq 18/19$ value in the table).
These figures emphasize the need to adapt the precision of the Krylov solver as a function of the current nonlinear residual.

Regarding local Newton iterations, we see that, except for the dual method which is inefficient in that case, no more than 3 times more computations of the tangent are required for the nonlinear approaches. 
In particular, the optimal Robin only requires 2.4 times more  computations of the tangent.
 \medskip

As a conclusion for this test case,  the primal and mixed versions of nonlinearly localized methods lead to interesting performance. Moreover they seem to be more effective when plasticity remains localized in only one subdomain. Indeed, global and Krylov numbers of iterations present a strong decrease when comparing them to those of the classic method (ratios from 0.61 to 0.76 at step $u_D=5u_e$), which means a decrease in communications between subdomains. The cost of this improvement stays limited, with less than 3 times more factorizations at step 5$u_e$ (and even less than 2 times more for optimal $\impe$). This is interesting because, as previously said, local iterations are independent and do not require communications between processors, which results in a gain of time for large -scale solutions.

The dual approach does not provide significant results, since additional increments had to be used to achieve convergence. Thus, it will not be considered in the following part, where better performance is achieved with an adaptation of stopping criteria.

\subsubsection{Adaptation of stopping criteria}

In order to observe the influence of the tuning of $\varepsilon_K$ and $\varepsilon_{NL}$ with respect to the current convergence state, the two expressions of section~\ref{section:sync} are tested with a set of coefficients chosen for their performance (see Table~\ref{tab:plast:param}). As discussed in  section~\ref{sec:error}, we choose $\varepsilon_K=\varepsilon_{NL}$. The same sequence of load increments as that of the previous section is applied to the perforated plate. 

\begin{table}[ht]
\begin{center}
\renewcommand{\arraystretch}{1.2}
\begin{tabular}{|c|c|c|}
\hline
Choice 1 & Eq~\eqref{expr:forcterm2} & $\gamma=0.7$, $\alpha=1.5$, $\varepsilon_{K_0}=10^{-6}$ \\
\hline
Choice 2 & Eq~\eqref{expr:forcterm3} & $\varepsilon_{K_0}=10^{-4}$ \\
\hline
\end{tabular}
\end{center}
\caption{Chosen parameters for $\varepsilon_K=\varepsilon_{NL}$}
\label{tab:plast:param}
\end{table}

The relative performance for choice 1 is summed up in Table~\ref{res:global_adapt_1}, whereas the relative performance for choice 2 is given in Table~\ref{res:global_adapt_2}. Depending on the quantity, we compare the variants with adapted criteria (INexact in the tables) to the fixed criteria approaches, or the nonlinear approaches with the classical approach, both with adapted criteria.
Note that choice 2 was too aggressive for classic and primal nonlinear localized methods to converge, so that only the mixed approaches were studied.

The aim of adaption is to avoid useless inner iterations without modifying the global iterations. This absence of perturbation corresponds to ratios of $1$ in the ``Global Newton'' rows. This happens almost every time.

Regarding the number of Krylov iterations, the figures in Table~\ref{res:global_adapt_1} are not as good, but still of the same order as the ones in Table~\ref{res:global_fix}. Table~\ref{res:global_adapt_2} shows that mixed approaches can be really efficient with a good adaption of criteria.  This means that adaption is not in contradiction with nonlinear localization and shall be used to reduce the number of exchanges.

Concerning local Newton iterations, which correspond to the extra cost of nonlinear localization techniques, we see that adaption makes it possible to reduce their number, with slightly smaller figures in Tables~\ref{res:global_adapt_1} than in Table~\ref{res:global_fix}, and more significant improvement in Table~\ref{res:global_adapt_2} for mixed approaches.

\begin{table}
\begin{center}
\renewcommand{\arraystretch}{1.5}
\begin{tabular}{|c|c|c|c|c|c|c|}
\hline
&Increment of loading & {$u_e$} &{2$u_e$} & {5$u_e$} & {5.75$u_e$} & {6.5$u_e$} \\
\hline
&Spread of nonlinearity & {Elastic} &  \multicolumn{2}{c|}{SD 1 } & SD 1-2 & SD 1-3\\
\hline
\multirow{4}{*}{\rotatebox[origin=c]{90}{Global Newton}}
&Classic IN/Classic & 1 & 1 & 1 & 1.06 & 1 \\
&Primal IN/Primal & 1 & 1 & 1 & 1 & 1 \\
&Mixed IN/Mixed, $\impe=\stiff_{bb}$ & 1 & 1 & 1 & 1 & 1.06 \\
&Mixed IN/Mixed, $\impe$ opti & 1 & 1 & 1 & 1.09 & 1.06 \\
\hline
\multirow{4}{*}{\rotatebox[origin=c]{90}{Krylov}}
&Classic IN/Classic & 1 & 1 & 0.77 & 0.77 & 0.69 \\
&Primal IN/Classic IN & 1 & 1 & 0.86 & 0.85 & 0.89 \\
&Mixed IN/Classic IN, $\impe=\stiff_{bb}$ & 1 & 1 & 0.86 & 0.85 & 0.91 \\
&Mixed IN/Classic IN, $\impe$ opti & 1 & 1 & 0.73 & 0.73 & 0.79 \\
\hline
\multirow{4}{*}{\rotatebox[origin=c]{90}{Local Newton}}
&Classic IN/Classic & 1 & 1 & 1 & 1.06 & 1 \\
&Primal IN/Classic IN & 1 & 2 & 2.54 & 2.56 & 2.78 \\
&Mixed IN/Classic IN, $\impe=\stiff_{bb}$ & 1 & 2 & 2.31 & 2.33 & 2.57 \\
&Mixed IN/Classic IN, $\impe$ opti & 1 & 2 & 1.85 & 2.06 & 2.30 \\\hline
\end{tabular}
\end{center}
\caption{Performance of adapted criteria with choice 1.}
\label{res:global_adapt_1}
\end{table}

\begin{table}
\begin{center}
\renewcommand{\arraystretch}{1.5}
\begin{tabular}{|c|c|c|c|c|c|c|}
\hline
&Increment of loading & {$u_e$} &{2$u_e$} & {5$u_e$} & {5.75$u_e$} & {6.5$u_e$} \\
\hline
&Spread of nonlinearity & {Elastic} &  \multicolumn{2}{c|}{SD 1 } & SD 1-2 & SD 1-3\\
\hline
\multirow{2}{*}{\rotatebox[origin=c]{90}{Global} \rotatebox[origin=c]{90}{Newton}}
&Mixed IN/Mixed, $\impe=\stiff_{bb}$ & 1 & 1 & 1 & 1 & 1.06 \\
&Mixed IN/Mixed, $\impe$ opti & 1 & 1 & 1 & 1 & 1 \\
\hline
\multirow{2}{*}{\rotatebox[origin=c]{90}{Krylov}}
&Mixed IN/Classic IN, $\impe=\stiff_{bb}$ & 1 & 0.79 & 0.71 & 0.67 & 0.66 \\
&Mixed IN/Classic IN, $\impe$ opti & 1 & 0.79 & 0.58 & 0.56 & 0.57 \\ \hline
\multirow{2}{*}{\rotatebox[origin=c]{90}{Local} \rotatebox[origin=c]{90}{Newton}}
&Mixed IN/Classic IN, $\impe=\stiff_{bb}$ & 1 & 2 & 2.23 & 2.28 & 2.43 \\
&Mixed IN/Classic IN, $\impe$ opti & 1 & 2 & 1.85 & 1.83 & 2.04 \\\hline
\end{tabular}
\end{center}
\caption{Performance of adapted criteria with choice 2.}
\label{res:global_adapt_2}
\end{table}

\subsection{Discussion}

Assessments proved that the nonlinear localization methods were able to reduce the number of Krylov iterations (and then the number of communication) at the price of some extra local computations (factorization), without modifying the number of global (outer Newton) iterations. Unfortunately our \texttt{octave-mpi} implementation does not allow reliable time measurements and large simulations. Moreover, the time performance of the methods would be highly dependent on the hardware and in particular on the performance of the network compared to the performance of the cores. 

We chose the previous experiment in order to prove that performance is highly dependent on the chosen boundary condition. In particular Neumann boundary conditions for perfectly plastic domains lead to very stiff problems, because in that case plasticity is estimated from the above (plasticity is overestimated and Newton iterations to reduce it are conducted on very unrealistic configurations). During our experimentations, we were able to design configurations with positive hardening where the dual approach was equivalent or even better than the primal. 

\begin{figure}[ht]
\centering
\includegraphics[width=.65\textwidth]{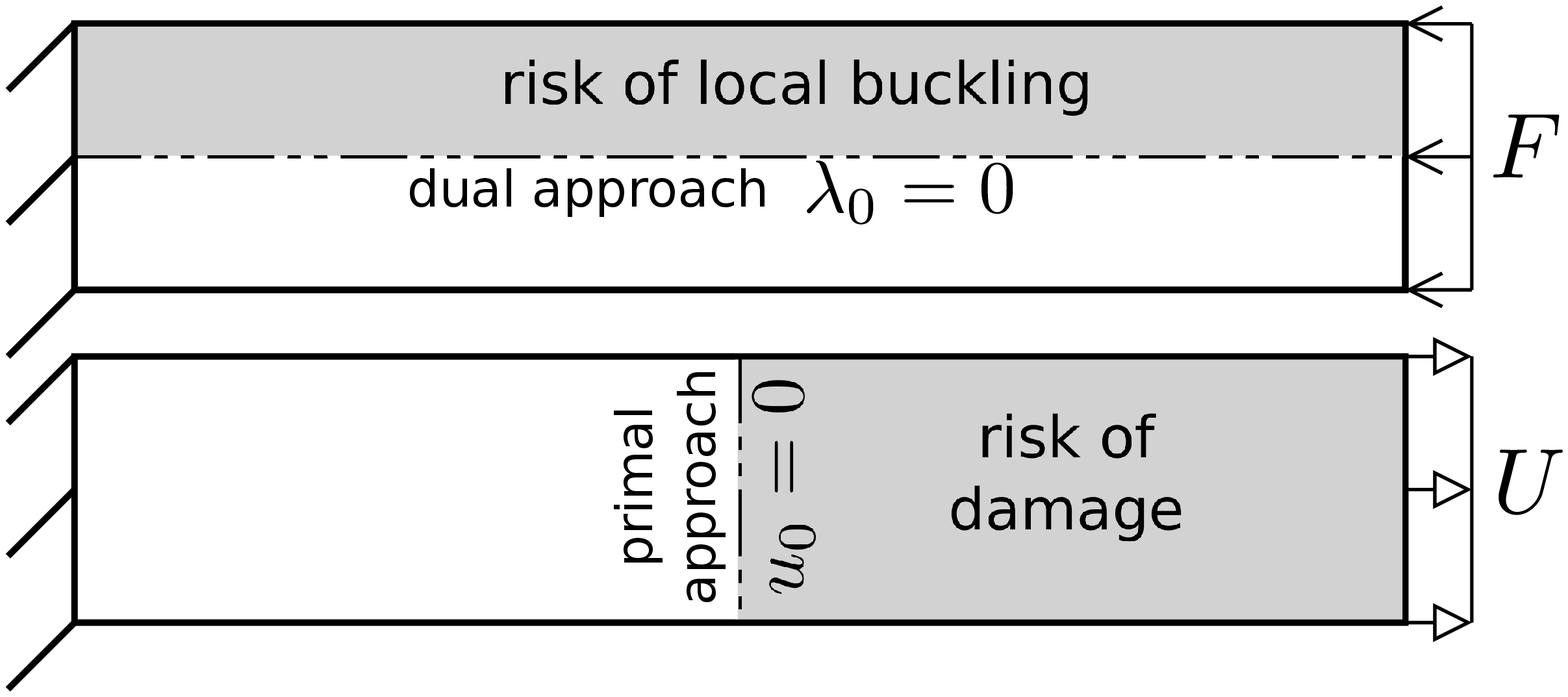}\caption{Two critical cases}\label{fig:simplefails}
\end{figure}

More generally, it is easy to build test cases where one variant behaves poorly, or even fails. Figure~\ref{fig:simplefails} presents two such cases. In the first case, the decomposition makes subdomains very slender and prone to buckling (assuming large transformations). In that case, the method may converge (slowly) to a suboptimal solution \cite{hinojosa14}. In the second case, the decomposition amplifies the initial deformation, and assuming damageable behavior, subdomains may fail even if the global problem was well-posed.

On the contrary, any experiment where the chosen boundary condition is
barely influenced by the nonlinearity will lead to better parallel performance than a linear problem: an extreme case is an isostatic elastoplastic lattice solved with a dual approach.

We can thus draw the following conclusions:
\begin{itemize}
\item The method may not converge or follow an inadequate path if local instabilities are possible. Based on the literature of Schwarz methods for nonlinear problems \cite{badea91,Lad07}, we expect the method to converge unconditionally (assuming small enough increments) in the case of coercive maximal monotone operators (typically positive hardening in small strains) where the existence and interesting properties of nonlinear Schur complements can be established \cite{showalter1977}. 
\item The shape of the subdomains shall play an important role in the potential development of nonlinear effects.
\item The initial boundary condition is critical, zero may be a very bad guess which triggers unduly complex local nonlinear computations. It is worth considering  starting the method by a global elastic prediction.
\item Well chosen Robin conditions should help follow the good path (the one of the monolithic solution). In the example above, we could compute the tangent Schur complement of the remainder of the structure which is (first order) optimal, but it is of pure academic interest, since the associated numerical cost would not be acceptable for actual industrial simulations; approximations of that operator is a classical question \cite{japhet01,GENDRE.2011.1} which we will address in future works. The approximation studied in the example, constituted by the stiffness of the neighbors, is not flexible enough, and it is not expected to give good results when the remainder of the structure is very compliant.
\end{itemize}

\section{Conclusion}
In this article, we have proposed nonlinear versions of FETI, BDD and FETI2LM algorithms where nonlinear computations are carried out independently on the subdomains. The starting point is a ``nonlinear condensation'' of the problems which leads, after linearization, to the parallel evaluation of the residual by nonlinear problems with boundary conditions characteristic of the method and to the solution to a classical DD system using any classical linear solver amongst FETI, BDD or FETI2NL.

The method modifies the nonlinear system to be solved and introduces inner inexact solvers, which implies a careful tuning of the stopping criteria. We proposed to exploit a classical inexact Newton formula.

A first assessment proved the potential of the method with an interesting trade-off between Krylov iterations (which corresponds to network communications) and Local Newton iterations (associated with local CPU computations). 

However the method suffers from several drawbacks. First, the performance depends on the chosen boundary condition, which is a strong difference with the linear case. Second, the shape of the subdomains plays an even greater role than in the linear case, since it can trigger unphysical nonlinear effects.  Third, in the case of localized nonlinearity, the method makes it possible to concentrate computations on the nonlinear subdomains but linear subdomains are then inactive, leading to poor load balancing (but potential power saving when idling processors enter standby mode). The question of load balancing will be addressed in future work when a more capable implementation is available. Typically we believe that it will be possible to propose an efficient strategy by encompassing the nonlinear problem in a mesh/substructuring adaption loop when a guaranteed precision is required by the user \cite{VREY.2013.1.2}.

Much work is thus still needed before the method is reliable. We believe that well-chosen Robin conditions are an important ingredient to make the method stable and efficient.

\bibliography{gosselet,biblio_paper,Biblio}
\end{document}